\documentclass{amsart}

\usepackage[latin1]{inputenc}
\usepackage{latexsym}
\usepackage{amsfonts}
\usepackage{amssymb}
\usepackage{color}
\usepackage{varioref}
\usepackage{verbatim}
\usepackage{epsfig}
\usepackage{hyperref}

\DeclareMathOperator{\kar}{char}
\DeclareMathOperator{\Res}{Res}
\DeclareMathOperator{\Spec}{Spec}

\newcommand{\secemail}{%
\setlength{\unitlength}{1pt}
bothmer
\begin{picture}(0,1)
\put(0,0){m}
\put(-4.2,0){@}
\end{picture}
ath.uni-hannover.de}

\newcommand{\xlabel}[1]{{\label{#1}}}

\newcommand{\barF}{\overline{\FF}}

\newtheorem{Thm}{Theorem}
\newtheorem{Lem}[Thm]{Lemma}

\newtheorem{Prop}[Thm]{Proposition}

\newtheorem*{Conj*}{Conjecture}
\newtheorem*{Thm*}{Theorem}

\theoremstyle{definition}
\newtheorem{Rem}[Thm]{Remark}
\numberwithin{Thm}{section}

\newcommand{\ZZ}{{\mathbb Z}}
\newcommand{\FF}{{\mathbb F}}
\newcommand{\KK}{{\mathbb K}}

\newcommand{\PP}{{\mathbb P}}

\newcommand{\QQ}{{\mathbb Q}}

\begin{document}

\title{The Casas-Alvero conjecture for infinitely many degrees}

\author[von Bothmer]{Hans-Christian Graf von Bothmer$^1$}
\address{Institut f\"ur Mathematik\\
          Universit\"at Hannover\\
          Welfengarten 1\\
          30167 Hannnover, Germany}
\urladdr{ \href{http://www-ifm.math.uni-hannover.de/~bothmer}{www-ifm.math.uni-hannover.de/\textasciitilde bothmer}}
\email{\secemail}
\thanks{$^1$ Supported by the Schwerpunktprogramm ``Global Methods in Complex
        Geometry'' of the Deutsche Forschungs Gemeinschaft.}

\author[Labs]{Oliver Labs$^2$}
\address{Mathematik und Informatik\\
         Geb\"aude E2.4\\
         Universit\"at des Saarlandes\\
         66123 Saarbr\"ucken, Germany}
\urladdr{ \href{http://www.OliverLabs.net}{www.OliverLabs.net}}
\email{\setlength{\unitlength}{1pt}
Labs
\begin{picture}(0,1)
\put(0,0){m}
\put(-4.5,0){@}
\end{picture}
ath.uni-sb.de, \setlength{\unitlength}{1pt}
mail
\begin{picture}(0,1)
\put(0,0){O}
\put(-4.1,0){@}
\end{picture}
liverLabs.net}
\thanks{$^2$ Supported by the Radon Institute for Computational and Applied
Mathematics (RICAM, Linz), Austrian Academy of Sciences.}

\author[Schicho]{Josef Schicho$^3$}
\address{Radon Institut for Computational and Applied Mathematics\\
          Austrian Academy of Sciences\\
          4040 Linz, Austria}
\urladdr{ \href{http://www.ricam.oeaw.ac.at/research/symcomp}{www.ricam.oeaw.ac.at/research/symcomp}}
\thanks{$^3$ Supported by the Austrian Science Fund (FWF)
        in the frame of the projects ``Solving Algebraic Equations'',
        P18992-N18, and SFB 013, subproject 03.}
\email{\setlength{\unitlength}{1pt}
josef.schicho
\begin{picture}(0,1)
\put(0,0){o}
\put(-4.2,0){@}
\end{picture}
eaw.ac.at}

\author[van de Woestijne]{Christiaan van de Woestijne$^4$}
\address{Institut f\"ur Mathematik B\\
         Technische Universit\"at Graz\\
         8010 Graz, Austria}
\urladdr{ \href{http://www.opt.math.tugraz.at/~cvdwoest}{www.opt.math.tugraz.at/\textasciitilde cvdwoest}}
\email{\setlength{\unitlength}{1pt}
c.vandewoestijne
\begin{picture}(0,1)
\put(0,0){t}
\put(-4.2,0){@}
\end{picture}
ugraz.at}
\thanks{$^4$ Supported by the Special Semester on Gr\"obner Bases and Related
Topics 2006, organised by RICAM Linz (in cooperation with RISC Hagenberg).}

\begin{abstract} 
Over a field of characteristic zero, it is clear that a polynomial 
of the form $(X-\alpha)^d$ has a non-trivial common factor with
each of its $d-1$ first derivatives. 
The converse has been conjectured by Casas-Alvero. 
Up to now there have only been some computational verifications for small
degrees $d$. In this paper the conjecture is proved 
in the case where the degree of the polynomial is a power of a prime number, or
twice such a power.  

Moreover, for each positive characteristic $p$, we give an example of a monic
polynomial of degree $d>p$ which is not a $d$th power but which has a common
factor with each of its first $d-1$ derivatives.  This shows that the
assumption of characteristic zero is essential for the converse statement to
hold.
\end{abstract}

\maketitle

%\tableofcontents{}
% \listoffigures{}
% \listoftables{}

\setlength{\parindent}{0pt}
\setlength{\parskip}{1ex plus 0.5ex minus 0.2ex}

%%%%%%%%%%%%%%%%%%%%%%%%%%%%%%%%%%%%%

%%%%%%%%%%%%%%%%%%%%%%%%%%%%%%%%%%%%%%%%%%%%
\section{Introduction}
\xlabel{sec:intro}

Let $\KK$ be a field and let $\KK[X]$ be the ring of univariate polynomials
over $\KK$. For any polynomial $P\in\KK[X]$ and for any nonnegative integer
$i$, we denote by $P^{(i)}$ the $i$th derivative of $P$, and by $P_i$ the $i$th
Hasse derivative, which is $P^{(i)}$ divided by $i\,!$ in characteristic zero.

This paper is concerned with the following question posed by E.\ Casas-Alvero
in characteristic zero:

\begin{Conj*}[Casas-Alvero] \xlabel{conjLinPot}
  Let $P$ be a monic univariate polynomial of degree $d$ over a field
  $\KK$.  Then $\gcd(P,P_i)$ is nontrivial for $i=1,\ldots,d-1$ if and only if
  $P=(X-\alpha)^d$ for some $\alpha\in\KK$.
\end{Conj*}

Note that the implication from right to left is trivial. The truth of the other
implication depends on the characteristic of the base field $\KK$.  For $d\le7$
and assuming that $\kar\KK=0$, the conjecture was proved in
\cite{DTGVdthPowConj}, using Gr\"obner basis computations. Since then
the authors of \cite{DTGVdthPowConj} have settled the case of $d=8$ as well
(personal communication). It seems that no other cases are known.

In this paper, we prove the conjecture in characteristic $0$ for infinitely
many degrees $d$. More precisely, in Section \ref{sec:proof1} we show

\begin{Thm*} %\xlabel{main}
  Let $d$ be of the form $p^k$ or $2p^k$ for some prime number $p$.
  Then the Casas-Alvero Conjecture holds in characteristic $0$ for polynomials
  of degree $d$.
\end{Thm*}

Since we could not find a reference for the fact that the Casas-Alvero
Conjecture does not always extend to characteristic $p$, we include
explicit counter-examples of degree $d>p$ for each $p$ in Section
\ref{sec:counter-examples}. 

This work grew out of discussions in a meeting of the authors at RICAM in the
frame of the special semester on Gr\"obner bases. The connection with 
Gr\"obner bases is made clear in Section \ref{sec:computational-aspects}, where
we discuss the computational aspects of the problem.

\subsection*{Notations and definitions.} 
We will write 
$$ 
  P = a_0X^d + a_1X^{d-1} + \ldots + a_d 
$$ 
throughout the paper. With this notation, we have
$$ 
  P_i = \binom{d}{i} a_0X^{d-i} + \ldots + \binom{i+1}{i}
  a_{d-i+1}X + \binom{i}{i} a_{d-i}  
$$ 
for the $i$th Hasse derivative of $P$. In characteristic $0$ one has $P_i =
P^{(i)}/i\,!$.

As $\binom{i}{i}=1$, we see that none of the Hasse derivatives vanishes
identically for all polynomials, regardless of the characteristic of the base
field, whereas in characteristic $p$, the $n$th usual derivative of any
polynomial is identically zero for all $n\ge p$. Thus, the conditions of the
Conjecture are strengthened in positive characteristic by using the Hasse
derivatives instead of the usual ones, whereas in characteristic $0$ there is
no difference.

%%%%%%%%%%%%%%%%%%%%%%%%%%%%%%%%%%%%%%%%%%%%
\section{Mixing characteristics, or There and Back Again} \xlabel{sec:proof1}

We now establish some results which will finally lead to a proof of our Theorem. 
Although the conjecture in general does not hold in positive
characteristic, it turns out that, for degrees $d$ and a prime $p$ as in the
Theorem, it is true in characteristic $p$, and this fact is used in the proof.

First, we note that the condition $\gcd(P,P_i)\ne 1$ is equivalent to the
\emph{resultant} $\Res_X(P,P_i)$ being zero (when $a_0\ne 0$). Let us
consider the ``generic'' polynomial 
$$
  P=a_0X^d + \ldots + a_d
$$
as an element of the polynomial ring
$$
  \ZZ[a_0,\ldots,a_d][X]; 
$$
and let us assume that $\gcd(P,P_i)\ne 1$ for $1\le i\le d-1$. It follows that
the vector of coefficients $(a_0,\ldots,a_d)$ belongs to an algebraic variety,
namely the set of common zeros of the equations $\{\Res_X(P,P_i) \mid 1\le i\le
d-1\}$. Therefore, we will study the set of points with coordinates in $\KK$ on
this variety, for any field $\KK$. 

Before this, however, we apply some simplifications: it is enough to consider
the conjecture for \emph{monic} polynomials only, so we assume $a_0=1$; and,
because we assume that $P$ has a common factor with its $(d-1)$st derivative,
which is linear, we know that $P$ has a zero in the base field. If we translate
this zero to $0$, we have $a_d=0$, and moreover, the property of $P$ of having
a common factor with its derivatives is preserved under this translation. 

The defining property is also invariant under scaling of the variable $X$; from
this, it follows that the equation $\Res_X(P,P_i)$ is \emph{homogeneous} of
weighted degree $d(d-i)$, if we give weight $j$ to the variable $a_j$, for
$1\le j\le d-1$. (It is consistent to give $a_0$ weight $0$ and $a_d$ weight
$d$, as well.) 

Putting this all together, the object of interest is the weighted projective
scheme 
\begin{equation}
  X_d \subseteq \PP_{\ZZ}(1,2,\ldots,d-1)
\end{equation}
over $\ZZ$, defined by the homogeneous ideal 
$$
  I_d = \left<\Res_X(P,P_i) \mid i=1,\ldots,d-1\right> \subseteq 
  R_d = \ZZ[a_1,\ldots,a_{d-1}],
$$
where $a_j$ has weight $j$ for $1\le j \le d-1$. We will consider the set
$X_d(\KK)$ of $\KK$-rational points on $X_d$ for any field $\KK$. 

Under the simplifications given above, if $P$ is a power of a linear
polynomial, then we must have $P=X^d$ and $a_1=\ldots=a_{d-1} = 0$. But this
trivial rational point $(0,\ldots,0)$ is excluded from $X_d$, as we consider
$X_d$ to be projective. Therefore, we have:

\begin{Prop} 
  The Casas-Alvero Conjecture holds for polynomials of degree $d$ over a field
  $\KK$ if and only if $X_d(\KK)$ is empty.
\end{Prop}

The following result permits us to draw conclusions about the situation in
characteristic $0$ from results in characteristic $p$. For any field $\KK$, we
write $\overline{\KK}$ for an algebraic closure of $\KK$; also, we consider the
\emph{base extension} $X_d\times\Spec\KK$ of $X_d$ to a projective scheme over
$\Spec\KK$. The scheme $X_d\times\Spec\KK$ is empty if and only if $X_d$ has no
points over $\overline{\KK}$. 

\begin{Prop} \xlabel{prop:frompto0}
  Let $d\ge 1$ be an integer. If $X_d(\barF_\ell)$ is empty for some prime
  $\ell$, then the Casas-Alvero Conjecture holds, for degree $d$, in
  characteristic $0$ and in characteristic $p$ for all but finitely many primes
  $p$.
\end{Prop}

\begin{proof}
As $X_d$ is weighted projective over $\ZZ$, the structure morphism
$$
  \phi_d : X_d \rightarrow \Spec(\ZZ)
$$
is proper \cite[Theorem II.4.9]{Hartshorne:77}.
In particular, the image $\phi_d(X_d)$ is closed. The complement
$U$ of $\phi_d(X_d)$ is exactly the set of points in $\Spec(\ZZ)$ where the
fiber under $\phi_d$ is empty. By assumption, $U$ is non-empty, because it
contains the prime $\ell$. It follows that $U$ is dense in $\Spec(\ZZ)$, and
hence contains the generic point, as well as all but finitely many primes.
But the fiber of the generic point is $X_d\times\Spec\QQ$, while the fiber of a
prime $p$ is $X_d\times\Spec\FF_p$, so that all these fibers are empty by 
definition of $U$.

By base extension, one sees easily that if $X_d\times\Spec\QQ$ is empty, then
so is $X_d(\KK)$ for any field $\KK$ of characteristic $0$, and the analogon
holds for characteristic $p$.
\end{proof}

We now come to some concrete statements about the problem. 

\begin{Prop} \xlabel{prop:case12}
  The Casas-Alvero Conjecture holds over any field in degrees $1$ and $2$.
\end{Prop}

\begin{proof}
  The linear case is trivial. Suppose $P=X^2+a_1X$ shares a factor with its
  derivative $2X+a_1$; then it follows easily that $a_1=0$, whether $2$ is
  equal to $0$ or not.
\end{proof}

We use a number-theoretical lemma. For an integer $a$ and a prime number $p$,
let $v_p(a)$ be the number of factors $p$ in $a$; we have $v_p(0)=\infty$.

\begin{Lem}
  Let $d$ be a positive integer, let $p$ be a prime number dividing $d$, and
  let $i$ be an integer with $0\le i\le d$. If $v_p(i)< v_p(d)$, then the
  binomial coefficient $\binom{d}{i}$ is $0$ modulo $p$.
\end{Lem}

\begin{proof}
  It is an old result of Kummer \cite[p.\ 116]{Kummer}, which is easily proved,
  that the highest power of $p$ that divides $\binom{d}{i}$ is equal to the
  number of `carries' that occur when $i$ and $d-i$ are added in $p$-adic
  notation. If $v_p(d)=e$, then the $p$-adic expansion of $d$ ends in $e$
  zeros, whereas by our assumption, the expansions of $i$ and $d-i$ end in less
  than $e$ zeros. Therefore, if we add them, we must incur at least one carry,
  and the binomial coefficient will be divisible by $p$.
\end{proof} 

\begin{Prop} \xlabel{prop:prime}
  Let $d\ge 1$ be an integer. If $d$ is a power of some prime number $p$,
  then $X_d(\barF_p)$ is empty.
\end{Prop}

\begin{proof}
Assume that $P\in\barF_p[x]=x^d+\ldots+a_{d-1}x$ is a polynomial
having a common factor with all its Hasse derivatives up to order $d-1$.
Note that by the Lemma, $\binom{d}{i}=0$ in $\FF_p$, for $i=1,\ldots,d-1$. In
particular $\binom{d}{d-1}=0$, hence $P_{d-1}=a_1$. The existence
of a common factor with $P$ implies that $a_1=0$. But then
$P_{d-2}=a_2$, hence $a_2=0$, and so on. Hence
$a_1=\ldots=a_{d-1}=0$, and because the origin is not contained
in the weighted projective space, it follows $X_d(\barF_p)$ is empty.
\end{proof}

The argument of the Proposition can be generalised as follows.

\begin{Prop} \xlabel{prop:induc}
  Let $d\ge 1$ and $k\ge 0$ be integers. If $d=np^k$ for some prime $p$, and if
  $X_n(\barF_p)$ is empty, then $X_d(\barF_p)$ is empty.
\end{Prop}

\begin{proof}
  The proof of the previous Proposition shows that $a_1=\ldots=a_{p^k-1}=0$.
  Now we have $P_{d-p^k}= \binom{d}{d-p^k} X^{p^k} + a_{p^k}$, where this
  time, the leading coefficient does not necessarily vanish in $\FF_p$. 
  Continuing, we see that again in $P_{d-p^k-1}$, the leading coefficient 
  vanishes, as well as the coefficient of $a_{p^k}$, and we obtain
  $a_{p^k+1}=0$. This process eventually shows that $a_i=0$ unless $p^k\mid i$.
  We obtain
  $$
    P = X^d + a_{p^k} X^{d-p^k} + \ldots + a_{d-p^k} X^{p^k},
  $$
  which in characteristic $p$ is equal to $Q^{p^k}$ for some polynomial 
  $Q\in\barF_p[X]$ of degree~$n$, because the field $\barF_p$ is
  perfect. This polynomial $Q$ again must have a common factor with all
  its Hasse derivatives up to order $n-1$, which is impossible if no
  such polynomials (except the trivial one $X^n$) exist in degree $n$.
\end{proof}

We can now prove the main result of this paper.

\begin{Thm*} %\xlabel{main}
  Let $d$ be of the form $p^k$ or $2p^k$ for some prime number $p$.
  Then the Casas-Alvero Conjecture holds in characteristic $0$ for polynomials
  of degree $d$.
\end{Thm*}

\begin{proof}
  If $d$ is a prime power $p^k$, then by Proposition \ref{prop:prime}, we see
  that $X_d(\barF_p)$ is empty. If $d=2p^k$, then we first invoke
  Proposition \ref{prop:case12} to show that no nontrivial quadratic examples
  exist in characteristic $p$, and then use Proposition \ref{prop:induc} to
  prove that $X_d(\barF_p)$ is empty in this case as well.

  We can now finish the proof by using Proposition \ref{prop:frompto0}, which
  allows us to conclude that $X_d(\KK)$ is also empty for any field $\KK$
  of characteristic $0$.
\end{proof}

%%%%%%%%%%%%%%%%%%%%%%%%%%%%%%%%%%%%%%%%%%%%
\section{Counter-Examples in Positive Characteristic} 
\xlabel{sec:counter-examples}

For each prime field $\FF_p$, we construct a monic polynomial $P$ of degree
$d>p$ that violates the Casas-Alvero Conjecture.

\begin{Prop} \xlabel{counterex}
  Let $p$ be a prime number, let $P = X^{p+1}-X^p\in\FF_p[X]$, and let
  $d=\deg P=p+1$. Then $P$ is not a $d$th power, but it has a non-trivial
  common factor with its Hasse derivatives $P_i$, $i=1, 2, \ldots, d-1$.
\end{Prop}

\begin{proof}
We have
$$
  P=X^{p+1}-X^p = X^{p}(X-1) \in \FF_p[X].
$$
Thus $P$ is not a $d$th power, and it has a common factor $X$ with $P_i$ for
$i=1,2,\ldots, d~-~2$. Moreover,
$$
  P_{d-1}=dX-1 \equiv X-1 \mod p,
$$
and $X-1$ divides $P$.
\end{proof}

\begin{Rem} 
  If we fix $d$, and assume that the Casas-Alvero Conjecture is true for
  degree $d$ in characteristic $0$, then it follows from Proposition
  \ref{prop:frompto0} that the primes $p$ for which counter-examples to 
  the Casas-Alvero Conjecture exist over $\FF_p$ are bounded. For example, for
  $d=3$ the Conjecture is true for a field of any characteristic, except $2$.
  By Proposition \ref{prop:induc}, this implies that the Conjecture holds in
  characteristic $0$ for all degrees under $30$, except possibly
  $12$, $20$, $24$, and $28$.

  However, the bound on the bad primes for a given degree $d$ may be quite
  large. For example, considering all quadrinomials of the form
  $X^6+aX^4+X^3+bX^2$ that possibly violate the conjecture, we find the
  counter-example
  $$
    P = X^6+3144481702696843 X^4+X^3+ 2707944513497181 X^2
  $$
  in characteristic $7390044713023799$, even though the conjecture holds for
  $d=6$ over $\QQ$.
\end{Rem}

%%%%%%%%%%%%%%%%%%%%%%%%%%%%%%%%%%%%%%%%%%%%
\section{Computational aspects}
\xlabel{sec:computational-aspects}

As was already said in the Introduction, the Casas-Alvero Conjecture may
in principle be verified computationally for any degree $d$. One way to do this
is to let a computer algebra package compute the polynomials $\Res_X(P,P_i)$,
for $i=1,\ldots,d-1$, and then compute a Gr\"obner basis of the ideal
$I_d$ generated by these resultants in $\QQ[a_1,\ldots,a_{d-1}]$. The
Conjecture is true for degree $d$ in characteristic $0$ if and only if the
Gr\"obner basis for $I_d$ (in any term ordering) contains, for each
$i$, an element whose leading term is a power of $a_i$ (see \cite{AL} for more
on these concepts).

We have done several such computations, using the packages {\sc Singular}
\cite{Singular}, {\sc Magma} \cite{Magma}, {\sc Macaulay 2} \cite{Macaulay2},
and {\sc Maple} \cite{Maple}, all of which offer Gr\"obner basis
computations for ideals in multivariate polynomial rings. 

However, the cost of these computations becomes prohibitive already for small
degrees $d$. Taking the smallest open case, $d=12$, one will find that even
computing the resultants is impossible, let alone computing a Gr\"obner basis.
The main problem here is that these are polynomials in many variables, which
tend to have exponentially many nonzero terms, whereas a Gr\"obner basis may
again contain exponentially many such polynomials.

Another approach, used by the authors of \cite{DTGVdthPowConj}, is to take the
zeros of $P$ as parameters instead of the coefficients. Here one obtains an
explicit decomposition of $I_d$ into ideals that are ``more primary'' than
$I_d$, and one proceeds with these. Unfortunately, the number of these
components grows exponentially with $d$, and the largest case that was solved
using this technique is $d=8$ (personal communication).

It may be expected that the computations become easier if we reduce the ideal
$I_d$ modulo a prime, because many terms may become zero. This is in fact the
case, and allows quick verification of the Conjecture for some small degrees.

Pushing the argument of Proposition \ref{prop:prime} a little further, one can
even show that if $d=p^k$, then the polynomials $\Res_X(P,P_i)$, for
$i=1,\ldots,d-1$, \emph{already form a Gr\"obner basis} of the reduction of the
ideal $I_d$ modulo $p$, for the weighted degree reverse lexicographic monomial
order with $a_1 < a_2 < \ldots < a_{d-1}$. The leading monomial of
$\Res_X(P,P_i)$ under this order is $a_{d-i}^d$. This provides an alternative
proof that the projective variety of the ideal over $\barF_p$ is empty.

However, when trying to solve the case $d=12$ by reducing modulo $p$, one finds
that there are counterexamples in characteristic $p=2,3,5,7,11,13$, while for
$p=17$ and larger we face the same complexity problems as in characteristic
$0$.

Another way to obtain partial information is to fix several coefficients a
priori. For example, when we put $P = X^6 + aX^4 + X^3 + bX^2$, as in the last
section, we can readily compute a Gr\"obner basis of $I_6$, even when we take
$\ZZ$ as the base ring instead of $\QQ$. The ideal $I_6$ is here generated by
three inhomogeneous equations in two variables (two resultants are zero), and its
Gr\"obner basis contains the integer 
$$ 
  M=13^3 \cdot 19^7 \cdot 67^2 \cdot 20771^2 \cdot 21379 \cdot 23993^3 \cdot
    7783207 \cdot 40362599 \cdot 7390044713023799.  
$$
It follows that some quadrinomial of the cited form violates the Conjecture in
characteristic $p$ if and only if $p$ is a prime factor of $M$.

\end{document}